\documentclass [12pt]  {article} 
\usepackage{amssymb}
\def \FF{{\mathbb{F}}}

\def \UUU{{\cal U}}
\def \Aut{{\rm Aut}}

\begin{document}

\begin{center}
{\Large {\bf Metacyclic groups as automorphism groups 
of compact Riemann surfaces}}\\
\bigskip
\bigskip
{\sc Andreas Schweizer}\\
\bigskip
{\small {\rm Department of Mathematics,\\
Korea Advanced Institute of Science and Technology (KAIST),\\ 
Daejeon 305-701\\
South Korea\\
e-mail: schweizer@kaist.ac.kr}}
\end{center}
\begin{abstract}
\noindent
Let $X$ be a compact Riemann surface of genus $g\geq 2$, and 
let $G$ be a subgroup of $\Aut(X)$. We show that if the Sylow 
$2$-subgroups of $G$ are cyclic, then $|G|\leq 30(g-1)$. 
If all Sylow subgroups of $G$ are cyclic, then, with two 
exceptions, $|G|\leq 10(g-1)$. More generally, if $G$ is 
metacyclic, then, with one exception, $|G|\leq 12(g-1)$. 
Each of these bounds is attained for infinitely many values 
of $g$. 
\end{abstract}

\subsection*{1. Introduction}

\noindent
If $G$ is a finite group, we write $|G|$ for its order 
and $G'$ for the subgroup generated by all commutators. 
Moreover, $C_n$ denotes the cyclic group of order $n$,
and $D_n$ is the dihedral group with $|D_n|=2n$, although
some of the sources we use would write $D_{2n}$ for that
group. 
\par
We are interested in subgroups $G$ of $\Aut(X)$ where $X$ 
always is a compact Riemann surface of genus $g\ge 2$, and 
$\Aut(X)$ is its group of conformal (i.e. orientation-preserving) 
automorphisms. It is a classical result that then 
$|\Aut(X)|\leq 84(g-1)$. See for example [Bre, Example 3.17]. 
If $|\Aut(X)|=84(g-1)$,  then $\Aut(X)$ is called a Hurwitz group.
There are infinitely many $g$ for which this can happen (see for 
example [Mb, Theorem 8]).
\par
If one imposes additional conditions on $G\subseteq \Aut(X)$, 
one usually gets smaller upper bounds on $|G|$. Below, in 
Theorem 2.3 we recall in detail the known result for the
case where $G$ is a metabelian subgroup of $\Aut(X)$. This already 
shows that if one wants a smooth upper bound such that for 
infinitely many values of $g$ this bound can be attained, then 
one should allow finitely many values of $g$ for which there
are exceptional groups $G$ that surpass that bound.
\par
In this sense, the optimal bound for solvable subgroups of 
$\Aut(X)$ is $|G|\leq 48(g-1)$; see [Ch], [G1]. 
If $G\subseteq \Aut(X)$ is supersolvable, then $|G|\leq 18(g-1)$
with one exception for $g=2$; see [Z2] and corrections in [GMl] 
and [Z3]. 
And by [Z1] we have $|G|\leq 16(g-1)$ if $G$ is nilpotent.  
The bounds $|G|\leq 4g+4$ and $|G|\leq 4g+2$ for abelian
respectively cyclic subgroups of $\Aut(X)$ are classical.
See also [G1]. 
\par
This list is not exhaustive. However, there is a class of groups
for which only partial results are known, namely the metacyclic 
groups. See the beginning of Section 4 for some clarifications 
on the definition of metacyclic.
\par
The best-known metacyclic groups are the dihedral groups, and 
for these we have $|G|\leq 4g+4$ by [BCGG, Corollary 2.6].
But certain series of $G\subseteq \Aut(X)$ described in [BJ] 
(see Theorem 2.4 below) show that for some metacyclic groups 
$|G|$ can go up to at least $12(g-1)$.
\par
More generally, for certain types of groups containing a cyclic
subgroup of index $2$ the paper [Mi] determines the smallest
genus on which such a group can act. These formulas are also 
very useful in our treatment of metacyclic groups with quotient
$C_4$ or $C_6$. 
\par
For general groups $G$ that have a cyclic quotient that is not 
too small, [MZ1] determine the smallest genus on which $G$ can 
act, possibly reversing the orientation. So the smallest genus 
on which $G$ can act conformally might at worst be bigger.
\par
The non-abelian groups of order $pq$ where $p$ and $q$ are odd 
primes are obviously split metacyclic. [W] determines every 
genus on which such a group can act.
\par
In Sections 3 to 5 we determine the optimal upper bounds for $|G|$ 
in terms of $g$ for four classes of groups: groups with cyclic 
Sylow $2$-subgroups, metacyclic groups, $Z$-groups (i.e., groups
whose Sylow subgroups are all cyclic), and, as a special case of 
that, groups of square-free order. The main results are Theorems
3.2, 4.8, 4.9  and 5.5 and Corollary 5.6.
\par
In Section 6 we use the optimal bound for $Z$-groups to improve
the upper bound from [Sch] on the exponent of a solvable group
$G\subseteq \Aut(X)$. Actually, the current paper grew out of 
a remark by the very helpful referee of my paper [Sch], who
suggested that it should be possible to find the optimal bound
for $Z$-groups in $\Aut(X)$ by going through the possible 
signatures of the corresponding Fuchsian groups.
\\

\subsection*{2. Known results}

\noindent 
Let $X$ be a compact Riemann surface of genus $g\geq 2$.
The the universal covering of $X$ is the complex upper 
halfplane $\UUU$ (or equivalently the unit disk) and $X$
is the quotient of $\UUU$ by the group of deck 
transformations. This can be refined as follows:
\par
If $G\subseteq \Aut(X)$, then there exists a Fuchsian group
$\Gamma=\Gamma(h;m_1 ,m_2 ,\ldots, m_r)\subseteq \Aut(\UUU)$,
and a torsion-free, normal subgroup $K$ of $\Gamma$, called 
a surface kernel, such that $\UUU/K\cong X$, $\Gamma/K\cong G$
and $\UUU/\Gamma\cong X/G$. One says that $G\subseteq \Aut(X)$ 
is covered by $\Gamma$. Moreover, $h$ is the genus of $X/G$.
The integers $m_i$ are $\geq 2$ and called periods. 
\par
We refer to [Bre, Chapter 1, Section 3] for a survey and only 
list the facts that we really need. The most important one is
\\ \\
{\bf Theorem 2.1.} \it
If $G$ is covered by 
$\Gamma(h;m_1 ,m_2 ,\ldots, m_r)$, then 
$$|G|=\frac{2}{2h-2+\sum_{i=1}^r (1-\frac{1}{m_i})}(g-1).$$
\rm
\\ 
A quick application, that will be very convenient for us,
is the following.
\\ \\
{\bf Lemma 2.2.} \it
\begin{itemize}
\item[(a)] If $|G|>4(g-1)$, then $\Gamma$ is a Fuchsian group
with quotient genus $h=0$ and $r\in\{3,4\}$.
\item[(b)] If $|G|>8(g-1)$, then $\Gamma=\Gamma(0;m_1 ,m_2 ,m_3)$
except possibly for $\Gamma(0;2,2,2,3)$ which gives $|G|=12(g-1)$.
\end{itemize}
\rm

\noindent
{\bf Proof.} \rm 
(a) [Bre, Lemma 3.18 (a)].\\
(b) One easily checks that $\Gamma(0;2,2,2,4)$ gives 
$|G|=8(g-1)$ and $\Gamma(0;2,2,3,3)$ gives $|G|=6(g-1)$.
\hfill$\Box$
\\ \\
As an abstract group $\Gamma(h;m_1, m_2 ,\ldots, m_r)$ 
is generated by $2h+r-1$ elements with certain relations.
Actually, thanks to Lemma 2.2 we only have to do calculations 
with such groups when $h=0$. In this case things simplify to
$$\Gamma(0;m_1, m_2 ,\ldots, m_r)\cong\langle 
x_1,\ldots,x_r\ |\ x_1^{m_1}=\cdots =x_r^{m_r}=1, 
x_1x_2\cdots x_r=1\rangle.$$
From this one can easily read off $\Gamma/\Gamma'$. This is 
very important for the corresponding $G$ we are interested 
in, because $G/G'$ is a quotient of $\Gamma/\Gamma'$. Note
however that $G'$ is not necessarily a quotient of $\Gamma'$,
as $G/G'$ might be a proper quotient of $\Gamma/\Gamma'$.
\par
We now present a known theorem on metabelian automorphism 
groups. It is representative for all the results that are 
briefly mentioned in the Introduction. The second reason 
is that we need some of its details in the proof of 
Lemma 4.6. Practically all theorems of this form were 
proved using Fuchsian groups. 
\par
When we say that a group $G$ acts on a Riemann surface
$X$, we always mean that it acts faithfully, that is, 
only the neutral element of $G$ acts as identity on $X$.
\\ \\
{\bf Theorem 2.3.} [ChP], [G2], [G3] \it
Let $X$ be a compact Riemann surface of genus $g\geq 2$.
Let $G\subseteq \Aut(X)$ be a metabelian group. Then 
$$|G|\leq 16(g-1)$$
with the following three exceptions:
\par
a group of order $24$ acting on $g=2$ and covered 
by $\Gamma(0;2,4,6)$;
\par
a group of order $48$ acting on $g=3$ and covered 
by $\Gamma(0;3,3,4)$;
\par
a group of order $80$ acting on $g=5$ and covered 
by $\Gamma(0;2,5,5)$.
\\
Conversely, there is a metabelian group of order $16(g-1)$
acting on a Riemann surface of genus $g\geq 2$ if and only 
if $g=2$ or $g=k^2 \beta +1$ is odd where $k$, $\beta$ are 
positive integers such that $\beta$ divides $1+\alpha^2$
for some integer $\alpha$.
\rm
\\ \\
{\bf Proof.} \rm 
The upper bound and the exceptions for $g=3$ and $g=5$ 
are from [ChP]. The exception for $g=2$, which had been 
overlooked in [ChP], is mentioned at the end of [G2]. 
See also [Bre, Example 18.5].
\par
The first proof that the bound $16(g-1)$ is attained for 
infinitely many $g$ is also in [ChP]. However, half of their 
$g$ are even; so by the criterion with $g=k^2\beta +1$ being
odd, which is [G3, Theorem 1.1], the corresponding $G$ cannot 
be metabelian. Thus [G3, Theorem 1.1]  corrects and completes 
the description of the occuring $g$.
\par
A description of the metabelian groups of order $16(g-1)$ 
in terms of generators and relations is given in 
[G3, Theorem 1.2].
\hfill$\Box$
\\ \\
If we try to prove an analogue of Theorem 2.3 for groups $G$ 
that are a semidirect product of two cyclic groups, or for
groups $G$ of square-free order, or for groups $G$ that have 
a quotient $C_8$, then the following three known series of 
groups of automorphisms outline the smallest upper bound 
for $|G|$ that we can possibly hope for. This will become 
important in the next three sections.
\\ \\
{\bf Theorem 2.4.} [BJ, Theorem 1 and Section 3] \it
\begin{itemize}
\item[(a)] For each prime $p\equiv 1\ mod\ 3$ there exists 
a Riemann surface $X$ of genus $p+1$ with some group 
$C_{2p}\rtimes C_6 \cong G\subseteq \Aut(X)$. So $|G|=12(g-1)$. 
Moreover, $G$ is covered by $\Gamma(0;2,6,6)$.
\item[(b)] For each prime $p\equiv 1\ mod\ 5$ there exists 
a Riemann surface $X$ of genus $p+1$ with some group
$C_p \rtimes C_{10} \cong G\subseteq \Aut(X)$. So $|G|=10(g-1)$.
Moreover, $G$ is covered by $\Gamma(0;2,5,10)$.
\item[(c)] For each prime $p\equiv 1\ mod\ 8$ there exists 
a Riemann surface $X$ of genus $p+1$ with some group
$C_p \rtimes C_8 \cong G\subseteq \Aut(X)$. So $|G|=8(g-1)$.
Moreover, $G$ is covered by $\Gamma(0;2,8,8)$.
\end{itemize}
\rm

\noindent
The following three individual groups of automorphisms will 
occur as exceptions in many of the theorems below. Therefore 
we find it worthwhile to give somewhat detailed descriptions.
\\ \\
{\bf Example 2.5.}
By [N, Theorem 2] there exists exactly one Riemann surface
of genus $2$ that has an automorphism of order $8$, namely
$y^2 =x^5 -x$. This is also called the Bolza surface and is
known to have automorphism group $GL_2(\FF_3)$. Compare 
[Bre, Example 18.5]. The Sylow $2$-subgroup is quasidihedral;
$$G\cong\langle a,b\ |\ a^8 =b^2 =1,\ bab^{-1}=a^3\rangle.$$
Moreover, [Bre, Example 18.5] shows that this is the only 
group of order $16$ that can act on genus $2$, and that it 
is covered by $\Gamma(0;2,4,8)$.
\rm
\\ \\
{\bf Example 2.6.} 
As a special case of [W, Corollary 4.2], the smallest genus
on which $G\cong C_7 \rtimes C_3$, the unique non-abelian 
group of order $21$, can act is $3$. Now it is known that 
there are exactly two Riemann surfaces of genus $3$ with an
automorphism of order $7$. The first one is the hyperelliptic
$y^2 =x(x^7 -1)$ (compare [N, Theorem 1]), whose full 
automorphism group is $C_{14}$. The other one is the Klein 
quartic with projective equation $x^3 y +y^3 z +z^3 x=0$.
It has $\Aut(X)\cong PSL_2(\FF_7)$. Then $G\cong C_7 \rtimes C_3$
is the normalizer of a Sylow $7$-subgroup in $\Aut(X)$.
That $G$ then is covered by $\Gamma(0;3,3,7)$. Compare the 
proof of Proposition 4.2 for the last claim.
\rm
\\ \\
{\bf Example 2.7.} 
By [Bre, Example 18.5] the non-abelian group $C_3 \rtimes C_4$
can act on at least one Riemann surface of genus $2$. It then
is covered by $\Gamma(0;3,4,4)$.
\rm
\\

\subsection*{3. Groups with cyclic Sylow 2-subgroups}

\noindent
The following results will be helpful in the next two sections,
but they might also be interesting in their own right. We start
by recalling a known fact.
\\ \\
{\bf Theorem 3.1.} \it
Let $G$ be a finite group with a cyclic Sylow $2$-subgroup 
$P$. Then the elements of odd order form a normal subgroup 
$N$ and $G\cong N\rtimes P$. In particular, $G$ is solvable.
\rm
\\ \\
{\bf Proof.} \rm 
It is a well-known corollary to Burnside's Transfer Theorem
(compare [H, Theorem 14.3.1] or [R, Theorem 10.1.8]) that
if the Sylow $p$-subgroup $P$ for the smallest prime divisor
of $|G|$ is cyclic, then it has a normal complement $N$. 
This shows $G\cong N\rtimes P$. In our case, $N$ obviously 
is formed by all elements of odd order.
\par
Now the solvability of $G$ is equivalent to the solvability
of $N$, which is assured by the Theorem of Feit and Thompson
[FT]. At the same time this shows that to get the solvability 
of $G$ one cannot avoid using (or implicitly proving) this 
monumental theorem.
\hfill$\Box$
\\ \\
Now we show that for this class of groups there exists 
an analogue to Theorem 2.3 and the results mentioned 
in the Introduction.
\\ \\
{\bf Theorem 3.2.} \it
Let $X$ be a compact Riemann surface of genus $g\geq 2$. 
Let $G$ be a subgroup of $\Aut(X)$ such that the Sylow 
$2$-subgroups of $G$ are cyclic. Then 
$$|G|\leq 30(g-1).$$ 
If $G$ attains this bound then $|G|\equiv 2\ mod\ 4$. 
\par
Conversely, there are infinitely many $g$, necessarily 
with $g\equiv 6\ mod\ 10$, for which this bound is 
attained.
\rm
\\ \\
{\bf Proof.} \rm 
The four biggest possible orders, $84(g-1)$, $48(g-1)$,
$40(g-1)$ and $36(g-1)$ are all divisible by $4$ but do 
not allow a quotient of order $4$. Compare 
[Bre, Lemma 3.18] or [GMl, Table 4.1]. 
So $|G|\leq 30(g-1)$. 
\par
If this bound is attained, the corresponding Fuchsian 
group is $\Gamma=\Gamma(0;2,3,10)$, which also only allows 
an abelian quotient $C_2$. Hence $|G|\equiv 2\ mod\ 4$, 
which is equivalent to $g-1$ being odd. Moreover,
$\Gamma'=\Gamma(0;3,3,5)$ and $\Gamma''=\Gamma(0;5,5,5)$. 
So $G'/G''\cong C_3$. By [H, Theorem 9.4.2] $G''/G'''$ 
cannot also be cyclic, so $G''/G'''\cong C_5 \times C_5$.
Hence $25$ divides $|G|$. So, besides being odd, $g-1$ 
is also divisible by $5$, i.e., $g\equiv 6\ mod\ 10$.
\par
If in [Ch, Theorem 3.2] we fix $m=5$ and take for 
$n$ any odd number, we get a Riemann surface of genus 
$5\cdot n^{12}+1$ with the action of a group of order 
$|G|=2\cdot 3\cdot 5^2\cdot n^{12}$.
\hfill$\Box$
\\ \\
If $|G|$ is divisible by a higher power of $2$, we get 
a smaller bound, even under weaker conditions.
\\ \\
{\bf Theorem 3.3.} \it
Let $X$ be a compact Riemann surface of genus $g\geq 2$
and $G\subseteq \Aut(X)$. If $G$ has a quotient $C_8$, then
$|G|\leq 8(g-1)$. This bound is sharp by Theorem 2.4 (c).
\rm
\\ \\
{\bf Proof.} \rm 
If $|G|>4(g-1)$ then by Lemma 2.2
the corresponding Fuchsian group must be 
$\Gamma(0;m_1,m_2,m_3)$ or 
$\Gamma(0;m_1,m_2,m_3,m_4)$. Moreover, to have a quotient
$C_8$, at least two of the periods must be divisible by 
$8$. The biggest $|G|$ one can get under these conditions
is $|G|=8(g-1)$ from $\Gamma(0;2,8,8)$; everything 
else is smaller.
\hfill$\Box$
\\ \\
Compare also [MZ1, Theorem 1].
\\ \\
{\bf Lemma 3.4.} \it
Let $X$ be a compact Riemann surface of genus $g\geq 2$. 
If $G\subseteq \Aut(X)$ has a quotient $C_4$, then 
$|G|\leq 10(g-1)$ except for the following cases:
\par
$|G|=12(g-1)$ for $\Gamma(0;3,4,4)$;
\par
$|G|=16(g-1)$ for $\Gamma(0;2,4,8)$;
\par
$|G|=12(g-1)$ for $\Gamma(0;2,4,12)$;
\par
$|G|=\frac{32}{3}(g-1)$ for $\Gamma(0;2,4,16)$.
\rm
\\ \\
{\bf Proof.} \rm 
If $|G|>10(g-1)$, then by Lemma 2.2 the corresponding
$\Gamma$ must be a triangle group $\Gamma(0;m_1,m_2,m_3)$.
Moreover, to have a quotient $C_4$, two of the periods 
must be divisible by $4$. This also excludes the exception 
$\Gamma(0;2,2,2,3)$.
\par
Now $\Gamma(0;2,4,4)$ is not a triangle group, 
$\Gamma(0;3,4,4)$ gives $|G|=12(g-1)$, 
$\Gamma(0;3,4,8)$ gives $|G|=\frac{48}{7}(g-1)$, 
and $\Gamma(0;4,4,4)$ gives $|G|=8(g-1)$.
Likewise $\Gamma(0;2,4,8)$ gives $|G|=16(g-1)$,
$\Gamma(0;2,4,12)$ gives $|G|=12(g-1)$,
$\Gamma(0;2,4,16)$ gives $|G|=\frac{32}{3}(g-1)$,
and $\Gamma(0;2,4,20)$ gives $|G|=10(g-1)$.
Finally, $\Gamma(0;2,8,8)$ gives $|G|=8(g-1)$.
\hfill$\Box$
\\ \\
{\bf Corollary 3.5.} \it
Let $X$ be a compact Riemann surface of genus $g\geq 2$. 
If $G\subseteq \Aut(X)$ has cyclic Sylow $2$-subgroups and
$4$ divides $|G|$, then $|G|\leq 12(g-1)$. 
\rm
\\ \\
{\bf Proof.} \rm 
If $8$ divides $|G|$ this is a corollary to Theorem 3.3. 
And if $8$ does not divide $|G|$, this excludes the case 
$|G|=16(g-1)$ in Lemma 3.4.
\hfill$\Box$
\\ \\
{\bf Remark 3.6.} 
If the Sylow $p$-subgroups of $G$ for all {\it odd} primes 
$p$ are cyclic, in contrast to Theorem 3.1 this does not 
suffice to guarantee that $G$ is solvable. Not even if 
in addition the Sylow $2$-subgroup has a cyclic subgroup 
of index $2$. For example, the simple groups $PSL_2(\FF_p)$
have this property.
\par
Bounds like in Theorem 3.2 do also not hold for such groups.
If $p\equiv \pm 1\ mod\ 7$, then by [Mb, Theorem 8] 
$PSL_2(\FF_p)$ is a Hurwitz group. See [Sch] for more 
elaborations on this.
\rm
\\

\subsection*{4. Metacyclic groups}

\noindent
A finite group $G$ is called {\bf metacyclic} if it has a cyclic
normal subgroup $N$ such that $G/N$ is also cyclic. 
A {\bf split metacyclic} group is a group $G\cong C_m\rtimes C_n$.
\par
The group $C_3\rtimes C_4$ where the elements of order $4$ act
as inversion on $C_3$ has a subgroup $C_3 \times C_2 \cong C_6$
of index $2$, but the extension of $C_6$ by $C_4 /C_2$ is non-split.
This shows that a metacyclic group can have several metacyclic 
structures and that a split metacyclic group can also have 
non-split metacyclic structures.
A metacyclic group that is not split metacyclic is for example 
the quaternion group $Q_8$ of order $8$.
\par
If $N\triangleleft G$ and $G/N$ are both cyclic, this implies 
of course that $G'\subseteq N$ and hence that $G'$ is cyclic. 
Some books, e.g. [H, p.146] use the more restrictive definition 
that a finite group $G$ is metacyclic if $G'$ and $G/G'$ are 
cyclic. For example, the metacyclic groups $C_p \times C_p$ or 
$D_4$ or $Q_8$ would not be metacyclic in that restrictive sense.
\\ \\
{\bf Lemma 4.1.} \it 
Every (split) metacyclic group of even order also has a (split) 
metacyclic structure with a quotient of even order.
\rm
\\ \\
{\bf Proof.} \rm 
If $G\triangleright N\cong C_{2^e m}$ and $G/N\cong C_n$ with
$m,n$ both odd, then $C_m \triangleleft G$ and $G/C_m =H$ has 
a normal subgroup $C_{2^e}$ with quotient $C_n$. By the 
Schur-Zassenhaus Theorem [R, Theorem 9.1.2] we have 
$H\cong C_{2^e} \rtimes C_n$. Moreover, $|\Aut(C_{2^e})|=2^{e-1}$;
so $C_n$ can only act trivially on $C_{2^e}$. Thus 
$H\cong C_{2^e n}$.
\par
By a similar argument we have 
$C_{2^e m}\rtimes C_n \cong C_m \rtimes (C_{2^e} \rtimes C_n)
\cong C_m \rtimes C_{2^e n}$.
\hfill$\Box$
\\ \\
Contrary to what the proof might seem to suggest upon 
the first quick reading, we cannot always make the normal 
cyclic subgroup odd. Just think of $C_2 \times C_2$.
\par
We now embark on proving an analogue of Theorem 2.3 for 
metacyclic groups.
\\ \\
{\bf Proposition 4.2.} \it
Let $X$ be a compact Riemann surface of genus $g\geq 2$
and $G\subseteq \Aut(X)$ a metacyclic group of odd order.
Then 
$$|G|\leq 9(g-1),$$ 
except for $G\cong C_7\rtimes C_3$ acting on genus $3$
(see Example 2.6).
\rm
\\ \\
{\bf Proof.} \rm 
By [MZ2] the only bigger odd orders are $15(g-1)$ and 
$\frac{21}{2}(g-1)$. In both cases $G/G'\cong C_3$ and
$G'$ is a quotient of $\Gamma(0;5,5,5)$ respectively 
$\Gamma(0;7,7,7)$. As $G'$ must be cyclic, this leaves 
only the possibilities $|G|=15$ and $|G|=21$. But a group 
of order $15$ is cyclic, and hence it cannot act on 
a Riemann surface of genus $2$. 
\hfill$\Box$
\\ \\
The best known (non-abelian) metacyclic groups are of 
course the dihedral groups, or more generally groups 
containing a cyclic subgroup of index $2$.
\\ \\
{\bf Theorem 4.3.} \it
If $G\subseteq \Aut(X)$ contains a cyclic 
subgroup of odd order and index $2$, then 
$$|G|\leq 6g.$$ 
This bound is attained for any 
$g\equiv 5\ mod\ 6$ by the groups 
$$G\cong\langle a,b |\ a^{3g}=b^2 =1,\ b^{-1}ab=a^{g-1}\rangle.$$
\rm
\noindent
{\bf Proof.} \rm 
Obviously $G\cong C_n \rtimes C_2$ where $n$ is odd. The minimum
genus $g^*$ on which such a group can act was determined in
[Mi, Theorem 3.3]. It also depends on the action of $C_2$ on 
$C_n$. If we take $G$  with generators and relations as in our 
theorem, we get $g^*=\frac{n}{3}$. Some easy estimates show that 
for the other cases in [Mi, Theorem 3.3] one cannot get a smaller 
$g^*$.
\hfill$\Box$
\\ \\
{\bf Theorem 4.4.} \it
If $G\subseteq \Aut(X)$ contains 
a cyclic subgroup of index $2$, then 
$$|G|\leq 8g.$$ 
This bound is attained for any $g$($\geq 2$) by the groups 
$$G\cong\langle a,b |\ a^{4g}=b^2 =1,\ b^{-1}ab=a^{2g-1}\rangle.$$
\rm
\noindent 
{\bf Proof.} \rm 
For the group with generators and relations as in our theorem
[Mi, Theorem 3.3] gives $g$ as the minimum genus on which it
can act. However, [Mi] only investigates a special type of
non-split extensions of $C_n$ by $C_2$, and there are still 
other ones. So to establish the bound on all $G$ that contain
$C_n$ of index $2$, we instead use that $n\leq 4g+2$ by 
[N, Theorem 1]. But, given $g$, the Riemann surface $X$ of 
genus $g$ with $C_{4g+2}\subseteq \Aut(X)$ is uniquely determined
(still by [N, Theorem 1]), namely it is $y^2 =x(x^{2g+1}-1)$,
and actually $C_{4g+2}=\Aut(X)$. By [N, Theorem 2] the second 
biggest size of a cyclic group is $C_{4g}$, which furnishes 
our upper bound.
\hfill$\Box$
\\ \\
If $G$ has a bigger cyclic quotient, one also obtains a good 
bound, even without the assumption that $G$ is metacyclic.
\\ \\
{\bf Theorem 4.5.} \it
Let $X$ be a compact Riemann surface 
of genus $g\geq 2$ and $G\subseteq \Aut(X)$.
\begin{itemize}
\item[(a)] If $G$ has a quotient $C_{10}$, then $|G|\leq 10(g-1)$. 
This bound is attained by the groups in Theorem 2.4 (b).
\item[(b)] If $G$ has a quotient $C_{2p}$ with a prime $p\geq 7$, 
then $|G|\leq 7(g-1)$. 
\end{itemize}
\rm

\noindent
{\bf Proof.} \rm 
If $|G|>4(g-1)$, then by Lemma 2.2 $G$ must be covered by 
a Fuchsian group with quotient genus $h=0$ and $r\in\{3,4\}$. 
Moreover, two periods must be divisible by $p$ and two periods 
must be even. So the biggest $|G|$ we can get comes from 
$\Gamma(0;2,p,2p)$ and is $\frac{4p}{p-3}(g-1)$.
\hfill$\Box$
\\ \\
Theorem 4.5 can also be extracted from [MZ1, Theorem 1].
But for quotients $C_4$ and $C_6$ only partial results
are known. In both cases we establish the bound we want 
for metacyclic $G$, but without assuming a priori that 
the quotient $C_4$ or $C_6$ is the one from the metacyclic 
structure. 
\\ \\
{\bf Lemma 4.6.} \it
Let $G\subseteq \Aut(X)$ be metacyclic. If $G$ has 
a (not necessarily cyclic) normal subgroup $N$ with 
$G/N\cong C_4$, then $|G|\leq 12(g-1)$. 
\rm
\\ \\
{\bf Proof.} \rm 
By Lemma 3.4 we only have to exclude the possibility
$|G|=16(g-1)$ with Fuchsian group $\Gamma=\Gamma(0;2,4,8)$.
This group has $\Gamma/\Gamma'\cong C_2 \times C_4$, so the
quotient from the metacyclic structure can only be $C_2$ or
$C_4$. But $C_2$ is excluded by Theorem 4.4, unless $g=2$.
By [Bre, Example 18.5] there is only one $G$ of order $16$
for $g=2$. This must be the group in Example 2.5, which
obviously has no quotient $C_4$. So for the rest of the 
proof we can assume $g\geq 3$.
\par
Now we know that $G$ must contain a cyclic subgroup of 
order $4g-4$. By [N, Theorems 3, 4 and 5] this is only 
possible if $4g-4\leq 3g+3$, i.e., if $g\leq 7$. But 
by Theorem 2.3 there are no metabelian groups of order 
$16(g-1)$ for $g=7$, $6$, $4$. And if $g=5$, then $G$ 
cannot have an element of order $16$ by [N, Theorem 5]. 
Finally, by [Bro, Table 5] there are only two groups 
of order $32$ that can act faithfully on genus $3$,
and they are not metacyclic.
\hfill$\Box$
\\ \\
{\bf Lemma 4.7.} \it
Let $G\subseteq \Aut(X)$ be metacyclic. If $G$ has 
a (not necessarily cyclic) normal subgroup $N$ with
$G/N\cong C_6$, then either $|G|=12(g-1)$ with $G$ 
being covered by $\Gamma(0;2,6,6)$ or $|G|\leq 9(g-1)$. 
\rm
\\ \\
{\bf Proof.} \rm 
If $|G|>8(g-1)$, then $G$ must be covered by a triangle 
group $\Gamma(0;m_1 ,m_2 , m_3)$ (Lemma 2.2) and the quotient
$C_6$ implies that at least two periods must be even and at 
least two periods must be divisible by $3$.
\par
If two periods are divisible by $6$, we have 
$\Gamma(0;2,6,6)$ or $|G|\leq 8(g-1)$, as 
$\Gamma(0;2,6,12)$ gives already $8(g-1)$ 
and $\Gamma(0;3,6,6)$ gives $6(g-1)$. 
\par
Alternatively, the three periods can be divisible by
$2$, $3$ and $6$. As $\Gamma(0;2,6,9)$ gives $9(g-1)$
and $\Gamma(0;3,4,6)$ gives $8(g-1)$, there only remains 
the series 
$$\Gamma(0;2,3,6k)\ \ \hbox{\rm with}\ \  k\geq 2,\ \ 
\hbox{\rm which gives}\ \ |G|=12\frac{k}{k-1}(g-1).$$
We have to exclude this series. So let us assume that
$G$ is covered by $\Gamma=\Gamma(0;2,3,6k)$. Then 
$\Gamma/\Gamma'\cong C_6$, and hence $G/G'$ is a quotient
of $C_6$. In particular, $G$ does not have a quotient of 
order $4$. Since $G$ is metacyclic, this means that $4$ 
cannot divide $|G|$. This also implies that $k$ must be odd.
\par
Taking Theorem 4.3 into account, we must have a metacyclic
structure $C_n \cong N\triangleleft G$ and $G/N\cong C_6$.
Note that the formula for $|G|$ shows that $k$ divides $n$. 
Let $n=p_1^{e_1}\cdots p_s^{e_s}$ with (odd) primes 
$p_1 <p_2<\ldots <p_s$. If the subgroup $C_2$ acts as 
identity on the subgroup $C_{p_i^{e_i}}$ of $C_n$, one can
divide by the other cyclic subgroups of $C_n$ and get 
a quotient $C_{p_i^{e_i}} \rtimes C_6$ in which $C_2$ is central. 
So it has a quotient $C_{p_i^{e_i}} \rtimes C_3$ although the 
biggest quotient of $\Gamma(0;2,3,6k)$ of odd order can only 
be $C_3$. Thus $C_2$ must act as inversion on all $C_{p_i^{e_i}}$,
and hence on $C_n$. In other words, $G$ contains a dihedral
subgroup $D_n$ of order $2n$ and index $3$.
\par
By an analogous reasoning, none of the primes $p_i$ can be
congruent to $2$ modulo $3$. Otherwise $C_3\subseteq C_6$ 
would have to act as identity on $C_{p_i^{e_i}}$ and we would
get a quotient $D_{p_i^{e_i}}$ of $\Gamma(0;2,3,6k)$.
\par
If $p_1^2$ divides $n$ or if $n$ is prime, then by 
[Mi, Corollary 3.4] the smallest genus on which $D_n$ 
can act is $g^*=n-\frac{n}{p_1}\geq \frac{2n}{3}$. Thus 
$|G|\leq 9g$. On the other hand, 
$12(g-1)<12\frac{k}{k-1}(g-1)=|G|$. This is only possible 
for $g\leq 3$. So $|G|\leq 27$, and hence $n\leq 4$, i.e.,
$n=3$, leading to $k=3$. Thus $|G|=18(g-1)$, which is too 
big for $g=3$ and does not exist for $g=2$ 
[Bre, Example 18.5].
\par
If $p_1^2$ does not divide $n$ and $n$ is also not prime
(so $n\geq 3\cdot 7=21$), then by [Mi, Corollary 3.4] the
smallest genus on which $D_n$ can act is 
$g^*=n+1-\frac{n}{p_1}-p_1\geq n+1-\frac{n}{3}-3=\frac{2n}{3}-2$,
so $|G|\leq 9(g+2)$. On the other hand, $12(g-1)<|G|$. This
is only possible for $g\leq 9$. Thus $|G|\leq 99$, giving
the contradiction $n\leq 16$.
\hfill$\Box$
\\ \\
Now we put everything together.
\\ \\
{\bf Theorem 4.8.} \it
Let $X$ be a compact Riemann surface of genus $g\geq 2$. 
If $G\subseteq \Aut(X)$ is metacyclic, then, with the 
exception of the  group of order $16$ described in
Example 2.5, we have
$$|G|\leq 12(g-1).$$
This bound is attained by the split metacyclic groups 
in Theorem 2.4 (a). 
\rm
\\ \\
{\bf Proof.} \rm 
If $|G|$ is odd, see Proposition 4.2.
\par
If $|G|$ is even, then in $C_m\cong N\triangleleft G$ and 
$G/N\cong C_n$ we can assume that $n$ is even. If $n$ is 
divisible by an odd prime, we can apply Theorem 4.5 or 
Lemma 4.7. If $4$ divides $n$, we take Lemma 4.6. Finally, 
if $n=2$, then by Theorem 4.4 we have $|G|\leq 8g$, which 
is $\leq 12(g-1)$ except for $g=2$. By [Bre, Example 18.5]
the only possible exception for $g=2$ is the group of order
$16$.
\hfill$\Box$
\\ \\
Accola [A] and Maclachlan [Ml] have shown independently that 
for every $g\geq 2$ there is a surface $X$ of genus $g$ with
$|\Aut(X)|\geq 8(g+1)$ and that this {\it lower} bound is sharp
for infinitely many $g$.
\par
The referee has suggested that Theorem 4.4 above might play 
a similar role for metacyclic groups. This is indeed the case.
\\ \\
{\bf Theorem 4.9.} \it
Let $p\geq 17$ be a prime that is congruent to $2$, $8$ or
$14$ modulo $15$. Then no Riemann surface of genus $g=p+1$
can have a metacyclic group $G\subseteq \Aut(X)$ with 
$|G|>8g$.
\rm
\\ \\
{\bf Proof.} \rm 
By [BJ, Theorem 1 and Section 6], if $g=p+1$ with a prime
$p\geq 17$, there are at most $6$ possibilities for 
$G\subseteq \Aut(X)$ with $|G|\geq 8(g-1)$. 
\par
Of these, $|G|=12(g-1)$ and $|G|=10(g-1)$ cannot occur because 
$p\not\equiv 1\ mod\ 3$ resp. $p\not\equiv 1\ mod\ 5$, and
$|G|=8(g-1)$ is too small. The groups $G$ of orders $8(g+3)$ 
and $8(g+1)$ are not metacyclic. The first one because is has
a non-metacyclic quotient $S_4$. The other one is covered by
$\Gamma(0;2,4,2p+4)$, which shows that $G/G'$ must be a quotient 
of $C_2 \times C_2$; so if $G$ were metacyclic it would have 
a cyclic subgroup of order $4g+4$, which is impossible ([N]). 
The remaining group $G$ of order $8g$ is the one we described 
in Theorem 4.4. 
\hfill$\Box$
\\

\subsection*{5. Z-groups}

\noindent
A finite group whose Sylow subgroups are all cyclic is called
a {\bf Z-group}. It is not immediately obvious that such groups
are always split metacyclic.
\\ \\
{\bf Theorem 5.1. (Zassenhaus)} [H, Theorem 9.4.3], 
[R, Theorem 10.1.10] \it
A $Z$-group that is not cyclic can be written as a semidirect product
$$C_m\rtimes C_n$$
where $(m,n)=1$ and $m$ is odd. In particular, such a group is 
split metacyclic.
\rm
\\ \\
As $Z$-groups are metacyclic, the bound in Theorem 4.8 of course 
holds for them. But the examples we know that reach this bound,
namely the groups in Theorem 2.4 (a), are not $Z$-groups. This 
becomes clear from the original description in [BJ, Theorem 1],
where they are described as $G$ is a split extension of $C_p$ 
by $C_6 \times C_2$. 
\par
On the other hand, the groups in Theorem 2.4 (b) are $Z$-groups.
So we now start to work towards the bound $10(g-1)$.
\\ \\
{\bf Proposition 5.2.} 
Let $X$ be a compact Riemann surface of genus $g\geq 2$, and 
$G\subseteq \Aut(X)$. If $G$ is a $Z$-group and $8$ divides $|G|$,
then $|G|\leq 8(g-1)$. This bound is attained by the groups from 
Theorem 2.4 (c).
\rm
\\ \\
{\bf Proof.} \rm 
This is immediate from Theorem 3.3.
\hfill$\Box$
\\ \\
{\bf Proposition 5.3.} \it
If $G\subseteq \Aut(X)$ is a $Z$-group with 
$|G|\equiv 4\ mod\ 8$, then, with the exception 
of $C_3 \rtimes C_4$ acting on genus $2$, we have
$|G|\leq 10(g-1)$.
\rm
\\ \\
{\bf Proof.} \rm 
We still have to eliminate the triangle groups $\Gamma(0;3,4,4)$ 
and $\Gamma(0;2,4,12)$ from Lemma 3.4. In view of Theorem 4.5 and 
Lemma 4.7 we have $G\cong C_m \rtimes C_4$ with odd $m$.
\par
Let us first assume that $G$ is covered by $\Gamma(0;3,4,4)$.
If $m$ is divisible by a prime $p\geq 5$, then $G$ has a quotient
of order $4p$. But the biggest quotient of $\Gamma(0;3,4,4)$ of
order prime to $3$ is $C_4$. So $m$ must be a power of $3$. 
Since $\Aut(C_{3^e})\cong C_{2\cdot 3^{e-1}}$, the subgroup $C_2$ in
$G$ is central. Dividing it out, we get a quotient $D_{3^e}$.
In that group every element order divides either $3^e$ or $2$.
So $D_{3^e}$ being a quotient of $\Gamma(0;3,4,4)$ means that it
is generated by  two elements of order $2$ whose product has 
order $3$. But that only gives $D_3$. So $e=1$ and 
$G\cong C_3 \rtimes C_4$, with $g=2$ because of $|G|=12(g-1)$.
\par
Now assume that $G$ is covered by $\Gamma(0;2,4,12)$. Obviously
$G$ has a quotient $H\cong C_3 \rtimes C_4$, which cannot be 
cyclic by Lemma 4.7. Since $H$ has only one involution, an
element of order $4$ and an element of order $2$ together only 
generate a group of order $4$. So $H$ and hence $G$ cannot be
a quotient of $\Gamma(0;2,4,12)$.
\hfill$\Box$
\\ \\
{\bf Proposition 5.4.} \it
If $G\subseteq \Aut(X)$ is a $Z$-group with 
$|G|\equiv 2\ mod\ 4$, then $|G|\leq 10(g-1)$.
\rm
\\ \\
{\bf Proof.} \rm 
If $G$ has a quotient $C_{2p}$ with a prime $p\geq 5$, this
follows from Theorem 4.5. If $G$ has a quotient $C_6$, the
condition $|G|\equiv 2\ mod\ 4$ excludes the possibility
$|G|=12(g-1)$ in Lemma 4.7. There remains the case that 
$G$ contains a cyclic group of odd order and index $2$. 
Then Theorem 4.3 yields $|G|\leq 6g$, which can only be 
bigger than $10(g-1)$ if $g=2$. But if $g=2$, the conditions
$|G|\leq 12$ and $|G|\equiv 2\ mod\ 4$ also give $|G|\leq 10$. 
\hfill$\Box$
\\ \\
Putting Propositions 5.2, 5.3, 5.4 and 4.2 together, 
we get the main result of this section.
\\ \\
{\bf Theorem 5.5.} \it
Let $X$ be a compact Riemann surface of genus $g\geq 2$.
If $G\subseteq \Aut(X)$ is a $Z$-group, then, with the exception 
of $C_3 \rtimes C_4$ acting on genus $2$ and $C_7 \rtimes C_3$
acting on genus $3$ (see Examples 2.6 and 2.7), we have 
$$|G|\leq 10(g-1).$$
This bound is attained by the groups in Theorem 2.4 (b). 
\rm
\\ \\
As groups of square-free order must obviously be
$Z$-groups, we also immediately obtain
\\ \\
{\bf Corollary 5.6.} \it
Let $X$ be a compact Riemann surface of genus $g\geq 2$.
If $G\subseteq \Aut(X)$ is a group of square-free order, 
then, with the exception of $C_7 \rtimes C_3$ acting on 
genus $3$, we have 
$$|G|\leq 10(g-1).$$
This bound is attained by 
the groups in Theorem 2.4 (b).
\rm
\\

\subsection*{6. More on the exponent of G}

\noindent
The exponent $exp(G)$ of a finite group $G$ is the least
common multiple of all element orders.
\par
In [Sch, Theorem 4.4] we showed that for $G\subseteq \Aut(X)$
the optimal upper bound on $exp(G)$ is $42(g-1)$. For the case 
of solvable $G$ we got a smaller bound in [Sch, Proposition 5.3]. 
For that we needed an upper bound for the order of a $Z$-group 
$G$ in $\Aut(X)$ (compare [Sch, Proposition 5.1]). 
\par
Recall in this context the easy fact that 
$exp(G)=\prod exp(G_p)$ 
where the product is over the different primes that divide
$|G|$ and $G_p$ is a Sylow $p$-subgroup of $G$. In particular, 
$exp(G)=|G|$ if and only if $G$ is a $Z$-group.
\par
Now that we have determined the optimal bound for $Z$-groups, 
we can improve the bound on $exp(G)$ for solvable $G$. We 
start with two special cases.
\\ \\
{\bf Lemma 6.1.} \it 
Let $X$ be a compact Riemann surface of genus $g\geq 2$, 
and let $G$ be a solvable subgroup of $\Aut(X)$. 
\begin{itemize}
\item[(a)] If $|G|=30(g-1)$, then $exp(G)$ divides $6(g-1)$.
\item[(b)] If $|G|=21(g-1)$, then $exp(G)$ divides $3(g-1)$.
\end{itemize}
\rm

\noindent
{\bf Proof.} \rm 
If $|G|=30(g-1)$, then $G$ is covered by $\Gamma(0;2,3,10)$.
So $G/G'\cong C_2$ and $G'$ is covered by $\Gamma(0;3,3,5)$.
This implies $G'/G''\cong C_3$ and $G''$ is a quotient of
$\Gamma(0;5,5,5)$. As $G''/G'''$ cannot also be cyclic by
[H, Theorem 9.4.2], we necessarily have 
$G''/G'''\cong C_5 \times C_5$. So the Sylow $5$-subgroup 
of $G$ cannot be cyclic, and the claim follows.
\par
The proof for (b) is practically the same.
\hfill$\Box$
\\ \\
{\bf Proposition 6.2.} \it 
Let $G\subseteq \Aut(X)$. If $G$ is solvable and the genus 
$g$ of $X$ is bigger than $2$, then
$$exp(G)\leq 12(g-1).$$
\rm

\noindent
{\bf Proof.} \rm 
If $exp(G)=|G|$, then $G$ is a $Z$-group, and hence 
$exp(G)\leq 10(g-1)$ by Theorem 5.2, except for 
$G\cong C_7\rtimes C_3$ on genus $3$.
\par
If $exp(G)=\frac{1}{2}|G|$, then $|G|$ cannot be 
$48(g-1)$, $40(g-1)$ or $36(g-1)$ because by 
[Sch, Proposition 5.3] we have $exp(G)\leq 16(g-1)$.
Lemma 6.1 excludes the remaining two solvable 
possibilities from [GMl, Table 4.1] with $|G|>24(g-1)$.
\par
If $exp(G)=\frac{1}{3}|G|$, then the Sylow $2$-subgroups 
of $G$ are cyclic, and hence 
$exp(G)\leq \frac{1}{3}30(g-1)=10(g-1)$ by Theorem 3.2.
\par
If $exp(G)\leq\frac{1}{4}|G|$, then 
$exp(G)\leq \frac{1}{4}48(g-1)=12(g-1)$.
\hfill$\Box$
\\ \\
{\bf Remark 6.3.} 
We don't know whether the bound in Proposition 6.2 is sharp.
But more precisely, the proof shows that for solvable $G$ with 
$g\geq 3$ the only possibilities with $exp(G)>10(g-1)$ are
$C_7 \rtimes C_3$ on genus $3$ and perhaps 
$exp(G)=12(g-1)$ for $|G|=48(g-1)$ or $|G|=24(g-1)$.
\\ \\
To see a bit more, let 
$$rad(|G|)=\prod_{p\ divides\ |G|}p$$
be the product of the different prime divisors 
of $|G|$, each one taken only with multiplicity 
one. Obviously $rad(|G|)$ divides $exp(G)$ and 
$exp(G)$ divides $|G|$.
\\ \\
{\bf Theorem 6.4.} \it 
Let $X$ be a compact Riemann surface of genus $g\geq 2$ 
and let $G\subseteq \Aut(X)$ be solvable. Then, except for 
$G\cong C_7\rtimes C_3$ on $g=3$, we have
$$rad(|G|)\leq 10(g-1).$$
This bound is attained by the groups from Theorem 2.4 (b).
\rm
\\ \\
{\bf Proof.} \rm 
If $rad(|G|)=|G|$, then $|G|$ is square-free and we can invoke 
Corollary 5.6. So we only have to show that for solvable $G$
with $|G|>20(g-1)$ we have $rad(|G|)\leq 10(g-1)$. By 
[GMl, Table 4.1] for solvable groups the orders bigger than
$20(g-1)$ are $48(g-1)$, $40(g-1)$, $36(g-1)$, $30(g-1)$,
$24(g-1)$, and $21(g-1)$. For four of them it is obvious
that $rad(|G|)\leq 10(g-1)$; the remaining two were treated 
in Lemma 6.1.
\hfill$\Box$ 
\\ \\
{\bf Remark 6.5.} 
If $G\subseteq \Aut(X)$ is not solvable, we have
$$rad(|G|)\leq 42(g-1)$$
from [Sch, Theorem 4.1]. Also by [Sch, Theorem 4.1], this bound 
is attained if and only if $G$ is a Hurwitz group of exponent 
$42(g-1)$ and this number is square-free. For $g>3$ the 
first condition is equivalent to $G\cong PSL_2(\FF_p)$ with 
a prime $p\equiv \pm 1\ mod\ 7$ and $g=\frac{p^3 -p}{168}+1$. 
(See [Sch, Theorem 4.4]). The second condition then is equivalent 
to $p$ being congruent to $\pm 3$ modulo $8$ and $\frac{p-1}{2}$ 
and $\frac{p+1}{2}$ both being square-free. So the question 
whether this bound is attained infinitely often is equivalent
to the number-theoretic question whether there are infinitely
many primes $p$ that are congruent to $\pm 13$ or to $\pm 27$ 
modulo $56$ such that $\frac{p-1}{2}$ and $\frac{p+1}{2}$ are
both square-free.
\\ \\ \\
{\bf Acknowledgements.} I would like to thank the referee of the 
paper [Sch] for a helpful remark that triggered the writing of 
the current paper. Also many thanks to the referee of the current 
paper for the careful reading, a constructive report, and for 
suggesting that I also try to prove something like Theorem 4.9. 
\\

\subsection*{\hspace*{10.5em} References}
\begin{itemize}

\item[{[A]}] R.~Accola: On the number of automorphisms of
a closed Riemann surface, \it Trans. Amer. Math. Soc.
\bf 131 \rm (1968), 398-408 

\item[{[BJ]}] M.~Belolipetsky and G.~A.~Jones: Automorphism 
groups of Riemann surfaces of genus $p+1$, where $p$ is prime,
\it Glasgow Math. J. \bf 47 \rm (2005), 379-393

\item[{[Bre]}] T.~Breuer: \it Characters and Automorphism Groups 
of Compact Riemann Surfaces, \rm LMS Lecture Notes 280, Cambridge 
University Press, Cambridge, 2000

\item[{[Bro]}] S.~A.~Broughton: Classifying finite group actions 
on surfaces of low genus, \it J. Pure Appl. Algebra \bf 69 
\rm (1990), 233-270 

\item[{[BCGG]}] E.~Bujalance, F.J.~Cirre, J.M.~Gamboa and G.~Gromadzki:
On compact Riemann surfaces with dihedral groups of automorphisms
\it Math. Proc. Cam. Phil. Soc. \bf 134 \rm (2003), 465-477 

\item[{[Ch]}] B.~P.~Chetiya: On genuses of compact Riemann 
surfaces admitting solvable automorphism groups, 
\it Indian J. Pure Appl. Math. \bf 12 \rm (1981), 1312-1318 

\item[{[ChP]}] B.~P.~Chetiya and K.~Patra: On metabelian groups of 
automorphisms of compact Riemann surfaces, \it J. London Math. Soc.
\bf 33 \rm (1986), 467-472

\item[{[FT]}] W.~Feit and J.~G.~Thompson: Solvability of groups 
of odd order, \it Pacific J. Math. \bf 13 \rm (1963), 775-1029 

\item[{[G1]}] G.~Gromadzki: Maximal groups of automorphisms of
compact Riemann surfaces in various classes of finite groups,
\it Rev. Real Acad. Cienc. Exact. Fis. Natur. Madrid \bf 82 no. 2
\rm (1988), 267-275 

\item[{[G2]}] G.~Gromadzki: On soluble groups of automorphism of 
Riemann surfaces, \it Canad. Math. Bull. \bf 34 \rm (1991), 67-73

\item[{[G3]}] G.~Gromadzki: Metabelian groups acting on 
compact Riemann surfaces, \it Rev. Mat. Univ. Complut. Madrid
\bf 8 no. 2 \rm (1995), 293-305

\item[{[GMl]}] G.~Gromadzki and C.~Maclachlan: Supersoluble
groups of automorphisms of compact Riemann surfaces, 
\it Glasgow Math. J. \bf 31 \rm (1989), 321-327

\item[{[H]}] M.~Hall: \it The Theory of Groups, \rm Macmillan,
New York, 1959 

\item[{[Mb]}] A.~M.~Macbeath: Generators of the Linear Fractional 
Groups, in: 1969 Number Theory (\it Proc. Sympos. Pure Math. 
\rm vol. XII, Houston, Tex., 1967) pp. 14-32

\item[{[Ml]}] C.~Maclachlan: A bound for the number 
of automorphisms of a compact Riemann surface, 
\it J. London Math. Soc. \bf 44 \rm (1969), 265-272 

\item[{[MZ1]}] C.~May and J.~Zimmerman: The symmetric genus 
of metacyclic groups, \it Topology and its Applications 
\bf 66 \rm (1995), 101-115

\item[{[MZ2]}] C.~May and J.~Zimmerman: The symmetric genus 
of groups of odd order, \it Houston J. Math. \bf 34 no. 2 
\rm (2008), 319-338

\item[{[Mi]}] G.~Michael: Metacyclic groups of automorphisms 
of compact Riemann surfaces, \it Hiroshima Math. J. \bf 31 
\rm (2001), 117-132

\item[{[N]}] K.~Nakagawa: On the orders of automorphisms of 
a closed Riemann surface, \it Pacific J. Math. \bf 115 no. 2 
\rm (1984), 435-443 

\item[{[R]}] D.~J.~S.~Robinson: \it A Course in the Theory of Groups,
\rm Springer GTM 80, New York - Berlin, 1982

\item[{[Sch]}] A.~Schweizer: On the exponent of the automorphism 
group of a compact Riemann surface, \it Arch. Math. (Basel) \bf 107 
\rm (2016), 329-340 

\item[{[W]}] A.~Weaver: Genus spectra for split metacyclic 
groups, \it Glasgow Math. J. \bf 43 \rm (2001), 209-218 

\item[{[Z1]}] R.~Zomorrodian: Nilpotent automorphism groups of 
Riemann surfaces, \it Trans. Amer. Math. Soc. \bf 288 no. 1 
\rm (1985), 241-255

\item[{[Z2]}] R.~Zomorrodian: Bounds for the order of supersoluble
automorphism groups of Riemann surfaces, \it Proc. Amer. Math. Soc.
\bf 108 no. 3 \rm (1990), 587-600

\item[{[Z3]}] R.~Zomorrodian: On a theorem of supersoluble 
automorphism groups, \it Proc. Amer. Math. Soc.
\bf 131 no. 9 \rm (2003), 2711-2713

\end{itemize}

\end{document}